\documentclass{article}
\usepackage{amssymb}

\usepackage{amsmath}
\usepackage{graphicx}

\begin{document}

\begin{center}
\textbf{FUNDAMENTAL DOMAINS OF GAMMA AND ZETA FUNCTIONS} \\[0pt]
\ \\[0pt]
CABIRIA ANDREIAN CAZACU$^{(a)}$ and DORIN GHISA$^{(b)}$ \\[0pt]
\ \\[0pt]
{\small (a) \emph{Simion Stoilow} Institute of Mathematics of Romanian
Academy, P.O.Box 1-764 RO-014700, Bucharest, Romania} \\[0pt]
{\small E-mail address: cabiria.andreian@imar.ro} \\[0pt]
{\small (b) York University, Glendon College, 2275 Bayview Avenue, Toronto,
Canada, M4N 3M6} \\[0pt]
{\small E-mail address: dghisa@yorku.ca}
\end{center}

{\small \bigskip }





\begin{abstract}
{\small Branched covering Riemann surfaces }$(\mathbb{C},f)${\small \ are
studied, where }$f${\small \ is the Euler Gamma function and the Riemann
Zeta function. For both of them fundamental domains are found and the group
of covering transformations is revealed. In order to find fundamental
domains, pre-images of the real axis are taken and a thorough study of their
geometry is performed. The technique of simultaneous continuation,
introduced by the authors in previous papers, is used for this purpose.
Color visualization of the conformal mapping of the complex plane by these
functions is used for a better understanding of the theory. For the Riemann
Zeta function the outstanding question of the multiplicity of its zeros, as
well as of the zeros of its derivative is answered.}
\end{abstract}

\textbf{1. Introduction}\bigskip

It has been proved in [2] that every neighborhood of an isolated essential
singularity of an analytic function contains infinitely many non overlapping
fundamental domains. In fact this is true as well for essential
singularities which are limits of poles or of isolated essential
singularities [2-4,6]. The Euler Gamma function and the Riemann Zeta
function have $\infty $ as their unique essential singularity. For the Gamma
function, $\infty $ \ is a limit of poles, while for the Zeta function it is
an isolated essential singularity. It follows that for each one of these
functions the complex plane can be written as a disjoint union of sets whose
interiors are fundamental domains, i.e. domains which are mapped conformally
by the respective function onto the complex plane with a slit [1], page 98.
Since there is a great deal of arbitrary in the choice of the fundamental
domains we can try to use the pre-image of the real axis in order to find
such a disjoint union of sets. As we will see next, this is working very
well for the Gamma function, while for the Zeta function a supplementary
construction is needed. The method of fundamental domains allows one to
extract a lot of information about the function, in particular about its
zeros, as well as the zeros of its derivative and to reveal global mapping
properties of the function. We have shown in [5] that to every rational
function of degree $n$, a partition of the complex plane into $n$ sets can
be associated such that the interior of every one of them is a fundamental
domain of the function. For transcendental functions, an important fact is
their behavior in the neighborhood of essential singularities, as it appears
in [2]. Moreover, since the fundamental domains are the \emph{leafs} of the
corresponding covering Riemann surface, they can be used in the study of the
group of covering transformations of the respective surface.

\bigskip

\bigskip \textbf{2. Global Mapping Properties of the Euler Gamma Function}

We use the explicit representation of the Euler Gamma function as a
canonical product [1]:

\begin{equation}
\Gamma (z)=(e^{-\gamma z}/ z)\prod_{n=1}^{\infty }(1+z/n)^{-1}e^{z/n},
\end{equation}

where $\gamma $ is the Euler constant

\begin{equation}
\gamma =\lim_{n\rightarrow \infty }(1+\frac{1}{2}+...+\frac{1}{n}-\ln
n)\approx 0.57722
\end{equation}

It is obvious from this representation that $\Gamma $ has the set of simple
poles $A=\{$ $0,-1,-2,...\}$ and has no zero. The product converges
uniformly on compact subsets of $\mathbb{C} \setminus A$ and therefore $w=$ $%
\Gamma (z)$ is a meromorphic function in the complex plane $\mathbb{C}$%
.\bigskip

\textbf{Theorem 1. }\textit{The extended plane }$\widehat{\mathbb{C}}$%
\textit{\ can be written as the union }$\widehat{\mathbb{C}}$\textit{\ }$%
=\cup _{n=-\infty }^{\infty }\overline{\Omega }_{n}$\textit{\ where }$\Omega
_{n}$\textit{\ are unbounded simply connected domains such that }$\Gamma $%
\textit{\ maps conformally every }$\Omega _{n}$\textit{\ onto }$\widehat{%
\mathbb{C}}\setminus L_{n},$\textit{\ where }$L_{n}$\textit{\ are slits
alongside the real axis}$.$\textit{\ These domains accumulate to infinity
and only there. The mapping }$\Gamma :\overline{\Omega }_{n}->\widehat{%
\mathbb{C}}\setminus \{0\}$\textit{\ is surjective for every }$n.$\textit{%
\bigskip }

\textit{Proof}: The number $\Gamma (x)$ is real for every real $x$ and the
graph of the function $x\rightarrow \Gamma (x)$ has the lines $x=0,$ $x=-1,$ 
$x=-2,$ $\ldots$ as vertical asymptotes [8]. \bigskip

\begin{center}
\includegraphics[width=2truein,height=2truein]{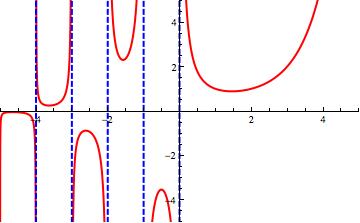} \\[0pt]
\ \\[0pt]
Figure 1
\end{center}

Fig. 1 can be found in most of the books of complex analysis serving as
texts for graduate studies. We used the online document [8]. It shows the
graph of the real function $x\rightarrow \Gamma (x),$ which can be used to
draw some information about the complex function $\Gamma .$

\bigskip

The respective graph has local minima and maxima, which correspond to the
points where $\Gamma ^{\prime }(z)=0.$ All these points are on the real
axis, namely $x_{0}\in (1,2),$ and for every positive integer $n,$ there is
a unique $x_{n}\in (-n,-n+1)$ such that $\Gamma ^{\prime }(x_{n})=0.$
Indeed, let us denote

\begin{equation}
\Gamma _{n}(z)=(e^{-\gamma z}/z)\prod_{k=1}^{n}(1+z/k)^{-1}e^{z/k}
\end{equation}

The sequence $(\Gamma _{n})$ converges uniformly on compact sets of $\mathbb{%
C}\setminus A$ to $\Gamma .$ It can be easily checked that for every $n\in
\mathbb{N}, $ the equation $\Gamma _{n}^{\prime }(x)=0$ is equivalent to an algebraic
equation of degree $n+1$ and has exactly $n+1$ real roots situated one in
every interval $(-k,-k+1),$ $k=1,2,...,n$ and one in the interval $(0,\infty
).$ Therefore $\Gamma _{n}^{\prime }(z)=0$ cannot have non real roots.

Since $(\Gamma _{n}^{\prime })$ converges in turn uniformly on compact
subsets of $\mathbb{C}\setminus A$ to $\Gamma ^{\prime },$ we infer that
every interval $(-n,-n+1),$ $n\in \mathbb{N},$ contains exactly one solution $x_{n}$
of the equation $\Gamma ^{\prime }(z)=0,$ and there is one\ more solution $%
x_{0}\in (1,2).$ There are no other solutions of this equation.

It is also obvious that:

\begin{equation}
\Gamma (x_{2k+1})<0\ \mbox{and} \ \Gamma (x_{2k})>0, k=0,1,2,\ldots
\end{equation}%
%

Based on this information, we can reveal the pre-image by $\Gamma $ of the
real axis, denoted $\Gamma ^{-1}(\mathbb{R})$. We'll see a little further
that all $x_{n},$ $n\geq 0$ are simple roots of $\Gamma ^{\prime }(z)=0,$
hence in a neighborhood $V_{n}$ of every $x_{n},$ $\Gamma (z)$ has the form
[1], page 133: 

\begin{equation}
\Gamma (z)=\Gamma (x_{n})+(z-x_{n})^{2}\varphi _{n}(z),\ \mbox{
where} \ \varphi _{n}(x_{n})\neq 0.
\end{equation}

By the Big Picard Theorem, the pre-image by $\Gamma $ of $\Gamma (x_{n})$ is
for every $n$ a countable set of points. The formula $(1)$ shows that $%
\Gamma (\overline{z})=\overline{\Gamma (z)},$thus this set is of the form $%
\{z_{n,k}\}\cup $ $\{\overline{z_{n,k}}\},k=0,1,2,...$ having the unique
accumulation point $\infty .$ Suppose that $z_{n,0}$ is $x_{n}.$ Then, the
pre-image of a small interval $(a_{n},b_{n})$ of the real axis centered at $%
\Gamma (x_{n})$ is the union of an interval $(\alpha _{n},\beta _{n})\ni
x_{n}$ of the real axis and another Jordan arc $\gamma _{-n}$ passing
through $x_{n},$ and symmetric with respect to the real axis, as well as
infinitely many other Jordan arcs passing each one through a $z_{n,k},$
respectively\ $\overline{z_{n,k}},$ $k>0.$ Simultaneous continuations [5]
over the real axis of these pre-images have as result the intervals $%
(-n,-n+1)$ for $(\alpha _{n},\beta _{n}),$ unbounded curves crossing the
real axis in $x_{n}$ for $\gamma _{-n}$ and infinitely many other unbounded
curves passing each one through $z_{n,k},$ or through $\overline{z_{n,k}}$
for $k>0.$ We use the same notation $\gamma _{-n}$ for the curves passing
through $x_{n}$ and $\gamma _{n,k}$, respectively $\overline{\gamma _{n,k}}$
for the others$.$ We notice that these curves cannot intersect each other,
since in such a point of intersection $z_{0}$ we would have $\Gamma ^{\prime
}(z_{0})=0,$ which is excluded. Also, the curves $\gamma _{n,k},$ and $%
\overline{\gamma _{n,k}}$ , $k\neq 0$ cannot intersect the real axis, for a
similar reason.\bigskip

\begin{center}
\includegraphics[width=2truein,height=2truein]{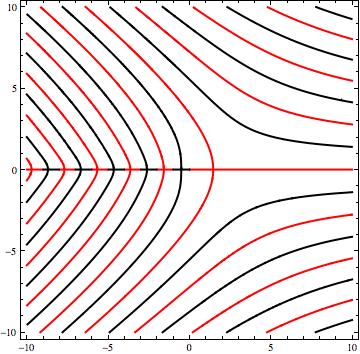} \\[0pt]
\ \\[0pt]
Figure 2
\end{center}

Fig. 2 shows the computer generated pre-image by $\Gamma $ of the real axis.
Since $0$ is a lacunary value for $\Gamma ,$ there can be no continuity
between the pre-image of the real positive half axis and negative half axis,
which means that each one of these pre-images has as components unbounded
curves \emph{closing} only at $\infty .$ Moreover, if we use two different
colors, say red and black, then these colors must alternate, since minima
and maxima for the real function $\Gamma (x)$ are alternating. Hence the
intervals $(\alpha _{n},\beta _{n})$ have alternating colors, which imply
alternating colors for $\gamma _{-n}.$Indeed, due to the continuity of $%
\Gamma ,$ except at poles, change of color can happen only at poles, which
means that the color of $\gamma _{-n}$ and that of $(\alpha _{n},\beta _{n})$
must agree. Thus, the colors of $\gamma _{-n}$ are alternating. On the other
hand, if a point travels on a small circle centered at origin in the $w$%
-plane ($w=\Gamma (z)$), it will meet alternatively the positive and the
negative half axis, which implies alternation into the colors of $\gamma
_{n,k}.$ If $x_{n}$ were multiple zeros of $\Gamma ^{\prime }$ then more
than one curve $\gamma _{-n}$ of the same color would start from $x_{n}$
violating this rule of color alternation. Thus, as previously stated, $%
\Gamma ^{\prime }$ has only simple zeros. The pre-image of the real axis
should be as we can see it in this computer generated picture.

\bigskip

We notice that if $z=x+iy\in \gamma _{n,k},$ or $z\in \overline{\gamma _{n,k}%
}$ then $\lim_{x\rightarrow -\infty }\Gamma (z)=0$ and $\lim_{x\rightarrow
+\infty }\Gamma (z)=\infty .$ Also, if $z\in \gamma _{-n},$ then $%
\lim_{x\rightarrow -\infty }\Gamma (z)=0.$ We will show next that the
domains bounded by some of the components of $\Gamma ^{-1}(\mathbb{R})$ are
fundamental domains of $\Gamma ,$ i.e. $\Gamma $ maps conformally each one
of them onto the complex plane with a slit. \bigskip

Let us first introduce notations for these domains. We denote by $G_{-n}$
the domains bounded by consecutive $\gamma _{-n},$ $n=0,1,2,...$ If $H_{+}$
and $H_{-}$ denote respectively the upper and the lower half plane, then
alternatively the image by $\Gamma $ of $H_{+}\cap G_{-n}$ and $H_{-}\cap
G_{-n}$ is the lower respectively the upper half plane. We denote by $\Omega
_{-n}=H_{+}\cap (G_{-2n-1}\cup G_{-2n}),$ $n\geq 0$. Let $\Omega _{1}$ be
the domain from the upper half plane bounded by $\gamma _{0},$ the interval $%
[x_{0},+\infty ]$ and the first component of $\Gamma ^{-1}(\mathbb{R})$
situated in the upper half plane and which does not intersect the real axis, 
$\Omega _{2}$ be the domain bounded by this component and the next one etc.
We denote by $\widetilde{\Omega }_{n},$ $n\in Z$ the domain symmetric to $%
\Omega _{n}$ with respect to the real axis.

\bigskip

Let us notice that the image by $\Gamma $ of every $\Omega _{-n}$ and of
every $\widetilde{\Omega }_{-n},$ $n\in N,$ is the complex plane with a slit
alongside the complement of respectively the interval $[\zeta (x_{2n-1},0]$
of the real axis, while the image of every $\Omega _{n}$ and of every $%
\widetilde{\Omega }_{n},$ $n=0,1,2,\ldots $ is the complex plane with a slit
alongside the positive real half axis. It is obvious that the domains $%
\Omega _{n}$ and $\widetilde{\Omega }_{n},$ $n\in Z$ accumulate to $\infty $
and only there in the sense that every neighborhood $V$ of $\infty $
contains infinitely many domains $\Omega _{n}$ and $\widetilde{\Omega }_{n}$
and any compact set in $\mathbb{C}$ intersects only a finite number of these
domains.

\bigskip

\begin{center}
\includegraphics[width=2truein,height=2truein]{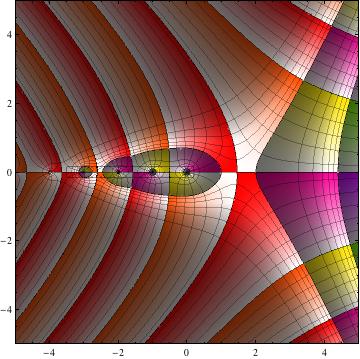} \\[0pt]
\ \\[0pt]
(a)\\[0pt]
\ \\[0pt]
\includegraphics[width=2truein,height=2truein]{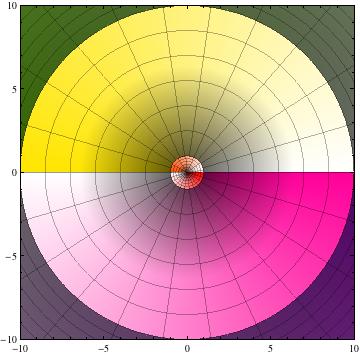}\ \ \ %
\includegraphics[width=2truein,height=2truein]{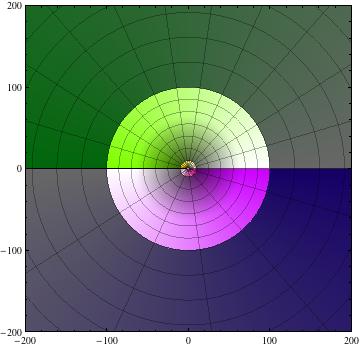}\ \ \ %
\includegraphics[width=2truein,height=2truein]{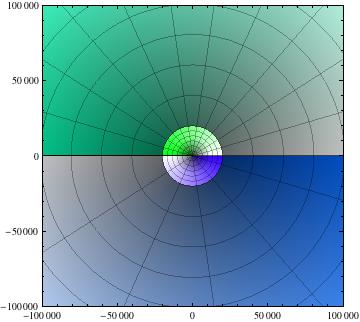}\\[0pt]
\ \\[0pt]
(b) \hspace{2in} (c)\hspace{2in} (d)\\[0pt]
\ \\[0pt]
Figure 3
\end{center}

Fig. 3 represents a visualization of the way the fundamental domains are
mapped conformally by $\Gamma $ onto the complex plane with a slit. Fig.
3(a) is obtained by taking pre-images of colored annuli centered at the
origin of the $w$-plane Fig. 3(b-d) and imposing the same color, saturation
and brightness on the pre-image of every point. The very big annuli Fig.
3(d) have pre-images around the poles and this is obvious when looking at
the colored pictures on the web project. However, the same colors appear for 
$z=x+iy$ with big positive values of $x$ characterizing the fact that $%
\lim_{x\rightarrow +\infty }\Gamma (x+iy)=\infty .$ Coupled with the
pre-image of orthogonal rays to these annuli, the picture Fig. 3(a) gives a
pretty accurate \emph{graphic} of the function.

\bigskip

\textbf{3. The Group of Covering Transformations of }$(\mathbb{C},\Gamma )$

\bigskip

Let us define the mappings $U_{k}:\mathbb{C}\rightarrow \widehat{\mathbb{C}}$
by setting

\begin{equation}
U_{k}(z)=\Gamma _{|\Omega _{k+j}}^{-1}\circ \Gamma (z)
\end{equation}

for every $z\in \Omega _{j},$ $k,j\in Z$ and by extending these mappings by
continuity to $\partial \Omega _{j}.$ In particular $U_{k}(-n)$ $=k-n,$ for $%
k-n\geqslant 0$ and $U_{k}(-n)=\infty $ for $k-n<0.$ Next, we extend $U_{k}$
to the lower half plane by symmetry:

\begin{equation}
U_{k}(\overline{z})=\overline{U_{k}(z)}
\end{equation}

We notice that $U_{k}$ are conformal mappings except for the points $-n$ and
for every $k\in Z$ we have $\Gamma \circ U_{k}(z)=\Gamma (z),$ $z\in \mathbb{%
C}.$ Moreover, $U_{k}(\Omega _{j})=\Omega _{k+j},$ $U_{k}(\widetilde{\Omega }%
_{j})=\widetilde{\Omega }_{k+j}.$

Finally we define $H:\mathbb{C}\rightarrow \widehat{\mathbb{C}}$ by :


\begin{equation}
H(z)=\Gamma _{|\widetilde{\Omega }_{k}}^{-1}\circ \Gamma (z)\ \mbox{if} \
z\in \Omega _{k}\ \mbox{and} \ H(z)=\Gamma _{|\Omega _{k}}^{-1}\circ \Gamma
(z),\ \mbox{if}\ z\in \widetilde{\Omega }_{k}
\end{equation}


and extend $H$ by continuity to $\partial \Omega _{k}$ and to $\partial 
\widetilde{\Omega }_{k},$ $k\in Z.$ It can be easily seen that $H$ is an
involution. We notice also that

\begin{equation}
U_{k}\circ U_{j}=U_{j}\circ U_{k}=U_{k+j}, U_{k}^{-1}=U_{-k},\ \ \ \ \ \ \
k,j\in Z,
\end{equation}

It is an elementary exercise to show that the group generated by $U_{1}$ \
and $H$ is the group of covering transformations of $(\mathbb{C},\Gamma ).$

\bigskip

Let us examine the function $\Gamma _{a}(z)=\Gamma (1/(z-a)).$ It is an
analytic function in $\mathbb{C}\setminus B,$ where

\begin{equation}
B=\{z\in \mathbb{C}:z=a\ \vee \ z=a-1/n,\ n\in N\}.
\end{equation}

Indeed, $\Gamma _{a}$ has no singular point in $\mathbb{C}\setminus B,$
since the transformation $z\rightarrow 1/(z-a)$ transforms the set $\{
0,-1,-2,\ldots \}\cup \{\infty \}$ into $B.$ The point $z=a$ is an essential
singular point and the points $a-1/n$ are poles of $\Gamma _{a}.$\bigskip

The domains $\Omega _{n}$ and $\widetilde{\Omega }_{n},$ $n\in Z$ are
transformed by $z\rightarrow 1/(z-a)$ into disjoint domains contracting
themselves to $a$ as $n\rightarrow \infty .$

Each one of these domains is mapped conformally by $\Gamma _{a}$ onto the
complex plane with a slit and the mapping of the closure of each one of them
on $\widehat{\mathbb{C}}\setminus \{0\}$ is surjective$.$

For arbitrary complex numbers $a_{1,}$ $a_{2,}$ $...,a_{m},$ let us define
now the function $f(z)=R(\Gamma _{a_{1}}(z),...,\Gamma _{a_{m}}(z))$ where $%
R $ is an arbitrary rational function. Let us denote $E=%
\{a_{1},a_{2},...,a_{n}\}$.

The function $f$ is meromorphic in $\widehat{\mathbb{C}}\setminus E$ and has 
$E $ as essential singular set. We state without proof the following.

\bigskip

\textbf{Theorem 2. }\textit{The set }$\widehat{\mathbb{C}}\setminus E$%
\textit{\ can be written as a disjoint union }$\widehat{\mathbb{C}}\setminus
E$\textit{\ }$=\cup H_{n}$\textit{\ where the interior }$\Omega _{n}$\textit{%
\ of every }$H_{n}$\textit{\ is mapped conformally by }$f(z)$\textit{\ onto }%
$\widehat{\mathbb{C}}\setminus L_{n},$\textit{\ where }$L_{n}$\textit{\ is a
part of a cut }$L.$ \textit{Every compact subset of }$\widehat{\mathbb{C}}%
\setminus E$\textit{\ meets only a finite number of }$\Omega _{n}$\textit{\
and every neighborhood of each }$a_{k}$\textit{\ includes infinitely many }$%
\Omega _{n}.$

\bigskip

\begin{center}
\includegraphics[width=2truein,height=2truein]{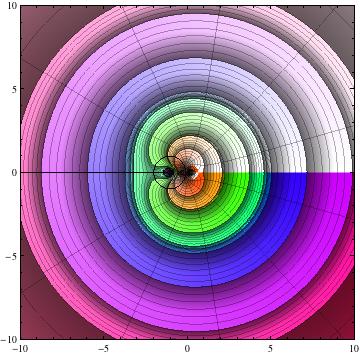} \ \ \ \ \ \ \ \ %
\includegraphics[width=2truein,height=2truein]{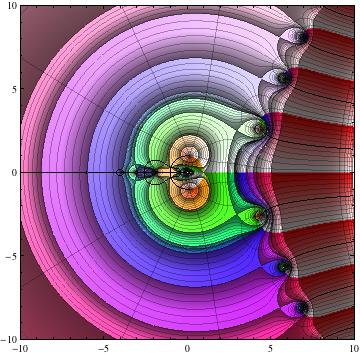}\\[0pt]
\ \\[0pt]
(a)\hspace{2in} (b) \\[0pt]
\ \\[0pt]
\includegraphics[width=2truein,height=2truein]{fig3c.jpg}\ \ \ \ \ \ \ \ %
\includegraphics[width=2.1truein,height=2truein]{fig3d.jpg}\ \ \\[0pt]
(c)\hspace{2in} (d)\\[0pt]
\ \\[0pt]
Figure 4
\end{center}

Fig. 4 represents a color visualization of the mapping realized by $\Gamma
(1/z)$ (Fig. 4(a)) and by $f(z)=\Gamma (z)+\Gamma (1/z)$ (Fig. 4(b)). Fig.
4(c,d) represents the annuli whose pre-images are shown in Fig. 4(a,b). Fig.
4 shows how the fundamental domains of $\Gamma (1/z)$ and those of $f(z)$
accumulate to the origin, which is an essential singularity for both of
them. We also can have an idea of the influence of one term of $f$ on the
graphic of the other term.

\textit{\bigskip }

\textbf{4. The Riemann Zeta Function}

\bigskip

Riemann Zeta function is one of the most studied transcendental functions,
in view of its many applications in number theory, algebra, complex
analysis, statistics, as well as in physics. Another reason why this
function has drawn so much attention is the celebrated Riemann conjecture
regarding its non trivial zeros, which resisted proof or disproof until now.

We are mainly concerned with the global mapping properties of Zeta function.
The Riemann conjecture prompted the study of at lest local mapping
properties in the neighborhood of non trivial zeros. There are known color
visualizations of the module, the real part and the imaginary part of Zeta
function at some of those points, however they do not offer an easy way to
visualize the general behavior of the function.\bigskip

The Riemann Zeta function has been obtained by analytic continuation [1],
page 178 of the series
\begin{equation}\zeta (s)=\sum_{n=1}^{\infty }n^{-s},\qquad s=\sigma
+it\end{equation}

which converges uniformly on the half plane $\sigma \geq \sigma _{0},$ where 
$\sigma _{0}>1$ is arbitrarily chosen. It is known [1], page 215, that
Riemann function $\zeta (s)$ is a meromorphic function in the complex plane
having a single simple pole at $s=1$ with the residue $1.$ Since it is a
transcendental function, $s=\infty $ must be an essential isolated
singularity. Consequently, the branched covering Riemann surface $(\mathbb{C}%
,\zeta )$ of $\mathbb{C}$ has infinitely many fundamental domains
accumulating at infinity and only there. The representation formula


\begin{equation}
\zeta (s)=-\frac{\Gamma (1-s)}{2\pi i} \int_{C}[(-z)^{s-1}/(e^{z}-1)]dz
\end{equation}


where $\Gamma $ is the Euler function and $C$ is an infinite curve turning
around the origin, which does not enclose any multiple of $2\pi i,$ allows
one to see that $\zeta (-2m)=0$ for every positive integer $m$ and there are
no other zeros of $\zeta $ on the real axis. However, the function $\zeta $
has infinitely many other zeros (so called, non trivial ones), which are all
situated in the (critical) strip $\{s=\sigma +it:0<\sigma <1\}.$ The famous
Riemann hypothesis says that these zeros are actually on the (critical) line 
$\sigma =1/2.$ Our study brings some new insight into this theory.

We will make reference to the Laurent expansion of $\zeta (s)$ for $|s-1|$ $%
>0:$


\begin{equation}
\zeta (s)=\frac{1}{s-1} +\sum\limits_{n=0}^{\infty }[(-1)^{n}/n!]\gamma
_{n}(s-1)^{n},
\end{equation}


where $\gamma _{n}$ are the Stieltjes constants:


\begin{equation}
\gamma _{n}=\lim_{m\rightarrow \infty }[\sum\limits_{k=1}^{m}(\log
k)^{n}/k-(\log m)^{n+1}/(m+1)]
\end{equation}


as well as to the functional equation [1], page 216:


\begin{equation}
\zeta (s)=2^{s}\pi ^{s-1}\sin \frac{\pi s}{2}\Gamma (1-s)\zeta (1-s).
\end{equation}

\bigskip

\textbf{5. The Pre-Image by }$\zeta $\textbf{\ of the Real Axis}

\bigskip

We will make use of the pre-image\ by $\zeta $\ of the real axis in order to
find fundamental domains for the branched covering Riemann surface $(\mathbb{%
C},\zeta )$ of $\widehat{\mathbb{C}}.$ By Big Picard Theorem, every value \ 
$z_{0}$ from the $z$-plane $(z=\zeta (s)),$ if it is not lacunary value, is
taken by the function $\zeta $ in infinitely many points $s_{n}$
accumulating to $\infty $ and only there. This is true, in particular, for $%
z_{0}=0.$

A small interval $I$ of the real axis containing $0$ will have as pre-image
by $\zeta $ the union of infinitely many Jordan arcs $\eta _{n,j}$ passing
each one through a zero $s_{n}$ of $\zeta ,$ and vice-versa, every zero $%
s_{n}$ belongs to some arcs $\eta _{n,j}.$ Since $\zeta (\sigma )\in \mathbb{%
R},$ for $\sigma \in \mathbb{R},$ and by $(15)$ the trivial zeros of $\zeta $
are simple zeros, the arcs corresponding to these zeros are intervals of the
real axis, if $I$ is small enough. We will show later that the non trivial
zeros of $\zeta $ are also simple, making superfluous the subscript $j$ in
the notation above. Due to the analiticity of $\zeta $ (except at $s=1$),
between two consecutive trivial zeros of $\zeta $ there is at least one zero
of the derivative $\zeta ^{\prime },$ i.e. at least one branch point of $(%
\mathbb{C},\zeta ).$ Thus, if we perform simultaneous continuation over the
real axis of the intervals $\eta _{n,j}$, we encounter at some moments these
branch points and the continuation follows on unbounded curves crossing the
real axis at these points.\bigskip\ 

Only the continuation of the interval containing the zero $s=-2$ stops at
the unique pole $s=1,$ since $\lim_{\sigma \nearrow 1}\zeta (\sigma )=\infty
.$ Similarly, if instead of $z_{0}=0$ we take another real $z_{0}$ greater
than $1$ and perform the same operations, since $\lim_{\sigma \searrow
1}\zeta (\sigma )=\infty ,$ the continuation over the interval $(1,\infty )$
stops again at $s=1.$ In particular, the pre-image by $\zeta $ of this
interval can contain no zero of the Zeta function. Thus, if we color red the
pre-image by $\zeta $ of the negative real half axis and let black the
pre-image of the positive real half axis, then all the components of the
pre-image of the interval $(1,+\infty )$ will be black,while those of the
interval $(-\infty ,1)$ will have a part red and another black, the junction
of two colors corresponding to a zero (trivial or not) of the function Zeta.

\bigskip

\begin{center}
\includegraphics[width=2truein,height=2truein]{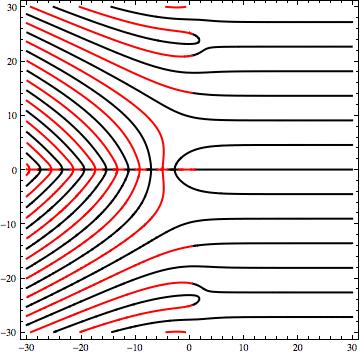} \\[0pt]
\ \\[0pt]
Figure 5
\end{center}

Fig. 5 represents the pre-image of the real axis in which the components
previously described are visible. We notice the existence of branch points
on the negative real half axis and their color alternation, as well as the
trivial zeros between them. Since these zeros are those of $\sin \frac{\pi s%
}{2},$ they are simple zeros and consequently there is no branching at them.
Some non trivial zeros are also visible. We'll explain next this
configuration.\bigskip

The red and the black unbounded curves passing through the branch points cannot
meet elsewhere (except $\infty $). Indeed, such an intersection point would
be a zero of $\zeta $ and the two curves would bound a domain which is
mapped conformally by $\zeta $ onto the complex plane with a slit alongside
the real axis from the image of one branch point to the image of the next
one. But such a domain should contain a pole of the function, which is
impossible.

The components of the pre-image of the real axis passing through non trivial
zeros form a more complex configuration. This configuration has something to
do with the special status of the value $z=1.$ Let us introduce notations
which will help making some order here and justifying the configurations
shown on the computer generated picture, Fig. 5. Due to the symmetry with
respect to the real axis, it is enough to deal with the upper half plane.
Let $x_{0}\in (1,+\infty )$ and let $s_{k}$ $\in \zeta ^{-1}(\{x_{0}\}).$
Continuation over $(1,+\infty )$ from $s_{k}$ is either an unbounded curve $%
\Gamma _{k}^{\prime }$ such that $\lim_{\sigma \rightarrow +\infty }\zeta
(\sigma +it)=1,$ by $(11),$ and $\lim_{\sigma \rightarrow -\infty }\zeta
(\sigma +it)=\infty ,$\ where $\sigma +it$ belongs to $\Gamma _{k}^{\prime
}, $ or there are points $u$ such that $\zeta (u)=1,$ thus the continuation
can take place over the whole real axis. We notice that it is legitimate to
let $\sigma $ tend to $-\infty $ on $\Gamma _{k}^{\prime }$, since if
infimum of $\sigma $ were a finite number $\sigma _{0},$ then $\sigma
_{0}+it_{0}\in \Gamma _{k}^{\prime }$ would be a pole of $\zeta ,$ which is
impossible.

Consecutive curves $\Gamma _{k}^{\prime }$ and $\Gamma _{k+1}^{\prime }$
form strips $S_{k}$ which are infinite in both directions. Indeed, if two
such curves met at a point $s,$ one of the domains bounded by them would be
mapped by $\zeta $ onto the complex plane with a slit alongside the real
axis from $1$ to $\zeta (s).$ Such a domain must contain a pole of $\zeta $,
which again cannot happen.

When the continuation can take place over the whole real axis, we obtain
unbounded curves containing each one a non trivial zero of $\zeta $ and a
point $u$ with $\zeta (u)=1.$ Such a point $u$ is necessarily interior to a
strip $S_{k}$ since the border of every $S_{k}$ and $\zeta ^{-1}(\{1\})$ are
disjoint. We denote by $u_{k,j}$ the points of $S_{k}$ for which $\zeta
(u_{k,j})=1$ and by $\Gamma _{k,j}$ the components of $\zeta ^{-1}(R)$
containing $u_{k,j}.$ The monodromy theorem assures that there is a one to
one correspondence between $u_{k,j}$ and $\Gamma _{k,j}.$

Let us notice that, when the continuation takes place over the whole real
axis, the components $\Gamma _{k,j}$ are such that the branches
corresponding to both the positive and the negative half axis contain only
points $\sigma +it$ with $\sigma <0$ for $|\sigma |$ big enough. Indeed, a
point traveling in the same direction on a circle $\gamma $ centered at the
origin of the $z$-plane meets consecutively the positive and the negative
real half axis. Thus the pre-image of $\gamma $ should meet consecutively
the branches corresponding to the pre-image of the positive and the negative
real half axis. On the other hand, due to the continuity of $\zeta $ on $%
\Gamma _{k}^{\prime },$ if a component of $\zeta ^{-1}(\gamma )$ meets a $%
\Gamma _{k}^{\prime },$ it should cross it and all $\Gamma _{l}^{\prime },$ $%
l\geq 1$ meeting consecutively the branches corresponding to the pre-image
of the positive and the negative half axis. Such an alternation is possible
only if the previously stated condition on $\sigma $ is fulfilled.

This analysis suggests that the value $z=1$ behaves simultaneously like a
lacunary value since $\lim_{\sigma ->\infty }\zeta (\sigma +it)=1,$ $\sigma
+it\in \Gamma _{k}^{\prime }$ and like an ordinary value, since $\zeta
(u_{k,j})=1.$ We can call it\textit{\ quasi-lacunary. }\bigskip

\bigskip

\begin{center}
\includegraphics[width=1.8truein,height=2truein]{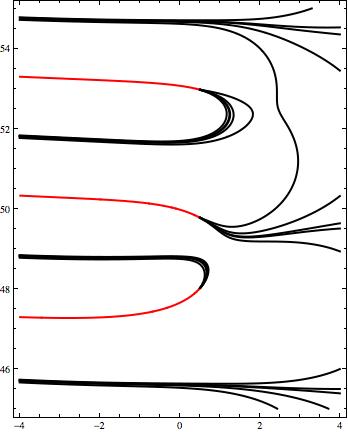} \\[0pt]
\ \\[0pt]
Figure 6
\end{center}

Fig. 6 represents dynamically the \emph{birth} of a strip. We picked up the
strip $S_{5}$ with $t$ in the range of $45$ and $55.$ It shows consecutively
domains which are mapped conformally by $\zeta $ onto the sectors centered
at the origin with angles from $\alpha $ to $2\pi -$ $\alpha ,$ where $%
\alpha $ takes respectively the values of $\pi /30,$ $\pi /100$ and $\pi
/1000.$ It is visible how the border of such a domain \emph{splits} into $%
\Gamma _{5}^{\prime },$ $\Gamma _{6}^{\prime }$ and $\Gamma _{5,0}$
previously defined as $\alpha \rightarrow 0.$\bigskip

It can be easily seen that two components $\Gamma _{k,j}$ cannot meet
neither, nor can they intersect any $\Gamma _{k}^{\prime }.$ Thus, those
components of pre-images of circles centered at the origin which cross a $%
\Gamma _{k}^{\prime },$ will continue to cross alternatively red and black
components of the pre-image of the real axis. These last components are
mapped by $\zeta $ either on the interval $(-\infty ,1),$ or on the whole
real axis.\bigskip

\textbf{Theorem 3. \ }\textit{For every }$k$\textit{\ there is a unique
component situated in the strip }$S_{k},$\ \textit{say }$\Gamma _{k,0}$%
\textit{, which is mapped bijectively by }$\zeta $\textit{\ onto }$(-\infty
,1).$\bigskip

\textit{Proof: }\ The strip $S_{k}$ is mapped by $\zeta $ onto the complex
plane with a slit alongside the real axis from $1$ to $+\infty .$ The
mapping is not necessarily bijective. For every $x_{0}\in (1,+\infty ),$
there is $s_{k}\in \Gamma _{k}^{\prime },$ and $s_{k+1}\in \Gamma
_{k+1}^{\prime }$ such that $\zeta (s_{k})=\zeta (s_{k+1})=x_{0}.$ Let us
connect $s_{k}$ and $s_{k+1}$ by a Jordan arc $\eta $ interior to $S_{k}$
(except for its ends)$.$ Then $\zeta (\eta )$ is a closed curve $C_{\eta }$
bounding a domain $D$ or a Jordan arc travelled twice in opposite
directions, in which case $D=\varnothing .$ We need to show that $C_{\eta }$
intersects again the real axis, in other words $\eta $ intersects the
pre-image of $(-\infty ,1)$. Indeed, otherwise $C_{\eta }$ would be
contained either in the upper, or in the lower half plane. Then $\zeta $
would map half of the strip $S_{k}$ bounded by $\eta $ and the branches of $%
\Gamma _{k}^{\prime }$ and $\Gamma _{k+1}^{\prime }$ corresponding to $%
\sigma \rightarrow +\infty $ onto $\mathbb{C}\setminus \overline{D}$ with a
slit alongside the real axis from $x_{0}$ to $1.$ We can take $x_{0}$ big
enough such that this half strip contains no zero of $\zeta ,$ which makes
impossible such a mapping.\bigskip

Let us show that $S_{k}$ cannot contain more than one component of the
pre-image of $(-\infty ,1).$ Indeed, if there were more, we could repeat the
previous construction with two consecutive such components, taking $s_{k}$
and $s_{k+1}$ with $\zeta (s_{k})=\zeta (s_{k+1})>0$ and arrive again to a
contradiction.

\bigskip

These arguments do not exclude the possibility of $S_{k}$ containing several
components $\Gamma _{k,j},$ $j\geq 1,$ which are mapped bijectively by $%
\zeta $ onto the whole real axis. Every one of these components contains a
non trivial zero of $\zeta $ and intersects the pre-image of the unit circle
in two points corresponding to $z=-1$ and $z=1.$ For $s=\sigma +it\in \Gamma
_{k,j}$ we have $\sigma \rightarrow -\infty $ for $\zeta (s)\rightarrow
\infty .$ If $S_{k}$ contains $m\geq 0$ such components we will call it $m$%
-strip. Every $m$-strip contains $m+1$ non trivial zeros of $\zeta $. Let us
denote by $s_{k,j}=\sigma _{k,j}+it_{k,j},$ $j=0,1,...,m$ \ the non trivial
zeros of $\zeta $ contained in $S_{k}.$ The computer generated data suggest
that $m$ is of the order of $\log t_{k,j}.$\bigskip\ \bigskip

\textbf{6. Fundamental Domains of the Riemann Zeta Function\bigskip }

The pre-image of circles centered at the origin of the $z$-plane are useful
in the study of the configuration of the components $\Gamma _{k,j}.$ The
circles with radius less than $1$ have bounded components of their
pre-images containing one or several zeros. The pre-image of the unit circle
has some bounded components containing each one a unique trivial or non
trivial zero and some unbounded components containing one or several such
zeros. All these curves must meet alternatively components of the pre-image
of the positive and negative real half axis. Indeed, a point moving in the
same direction on a circle centered at the origin will cross alternatively
the positive and the negative real half axis. A corollary of this fact is
that the branch points which are not zeros of $\zeta $ are simple zeros of $%
\zeta ^{\prime }.$This is due to the fact that in a neighborhood of such a
point all the components of the pre-image of the real axis must have the
same color and the rule of alternation would be violated if such a point
were multiple zero of $\zeta ^{\prime }.$ Thus, we have:\bigskip

\textbf{Theorem 4}.\textit{\ All the zeros of }$\zeta ^{\prime }$\textit{\
which are not zeros of }$\zeta $\textit{\ are simple zeros of }$\zeta
^{\prime }.$\bigskip

We will show later on that the zeros of $\zeta $ are all simple and
therefore all the zeros of $\zeta ^{\prime },$ with no exception, must be
also simple.

The pre-images of circles centered at the origin of radius less than or
equal to $1$ cannot meet the curves $\Gamma _{k}^{\prime }$ (which belong to
the pre-image of $(1,+\infty )$). We will see later that there are bounded
components of the pre-image of circles of radius greater than $1,$ but close
to $1$ with the same property. However, in the alternation of the components
of the pre-image of the positive and negative half axis for the other
circles, $\Gamma _{k}^{\prime }$ must be taken into account.\bigskip

\begin{center}
\includegraphics[width=1truein,height=2truein]{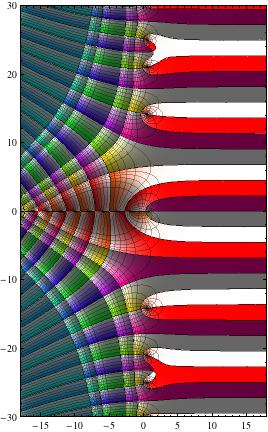}\ \ \ \ \ \ \ \ %
\includegraphics[width=2truein,height=2truein]{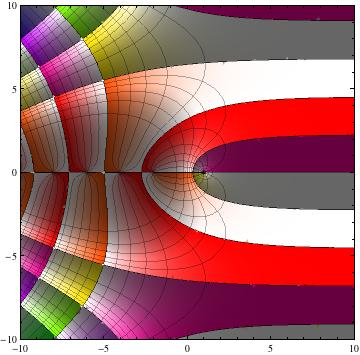}\hspace{1in} %
\includegraphics[width=.6truein,height=2truein]{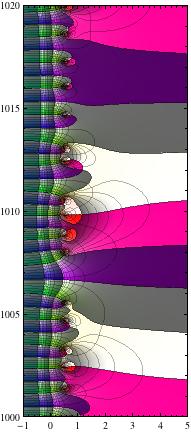}\\[0pt]
\ \\[0pt]
\hspace*{-5cm} (a) \hspace{1.6in} (b) \\[0pt]
\ \\[0pt]
\hspace{4cm} Figure 7 \hspace{2.25in} Figure 8
\end{center}

Fig. 7 shows the pre-images by $\zeta $ of the colored annuli from Fig. 3
intersecting the pre-image of the real axis in the box [-15,15] x [-30,30]
in Fig. 7(a) with a zoom on the origin in Fig. 7(b). The curves on the left
side in Fig. 7(a) crossing alternatively components of the pre-image of the
negative and positive real half axis are pre-images of circles centered at
the origin with radius greater than $1.$ The pre-image of annuli coupled
with the pre-image of some orthogonal rays give a pretty accurate
description of the mapping.\bigskip

Fig. 8 displays a $6$-strip situated in the area corresponding to $t\in
(1005,1016).$ There are clearly visible two components of the pre-image of
the unit circle: one bounded situated in the upper part of the strip
containing a unique non trivial zero, and one unbounded containing the other 
$6$ non trivial zeros. We notice in the strip above this $6$-strip two
bounded components of the pre-image of the unit circle. It appears that the
number of these bounded components also increases with $t.$\bigskip

Let $\gamma _{\rho }$ be a circle $|z|$ $=\rho $ for a small enough value of 
$\rho ,$ such that the pre-image of $\gamma _{\rho }$ is formed by disjoint
closed curves. If such a curve $\eta _{k,j}$ is in the critical strip, we
can suppose that it contains a unique non trivial zero $s_{k,j}$ of $\zeta .$
Suppose that more than one $s_{k,j}$ belong to the same component of the
pre-image of the unit circle. As $\rho $ \ increases the corresponding
curves $\eta _{k,j}$ expand such that for some value of $\rho ,$ $\eta
_{k,j} $ will meet another curve of the same type at a point $v_{k,j}.$
Indeed, all of them tend to the pre-image of the unit circle, as $\rho
\rightarrow 1.$ It is obvious that $v_{k,j}$ must be a branch point of $%
\zeta ,$ due to the fact that $\zeta $ takes the same value in points
situated on different curves $\eta _{k,j}$ in every neighborhood of $v_{k,j}$%
. Since $v_{k,j}$ cannot be a multiple pole, we have necessarily that $\zeta
^{\prime }(v_{k,j})=0.$ Let us notice that a curve $\eta _{k,j}$\ cannot 
\emph{split} into two or more such curves when $\rho $ suffers a small
change, since then two of them would meet in at least two different points
and would border a domain which is mapped by $\zeta $ onto the complex plane
with a slit alongside of an arc of $\gamma _{\rho },$ which is absurd.
Consequently, all the points $v_{k,j}$ must be simple zeros of $\zeta
^{\prime },$ otherwise we would have a \emph{split} of an $\eta _{k,j}$ at $%
v_{k,j}.$ Finally, since for an $s$ in the pre-image of the unit disc with $%
\zeta ^{\prime }(s)=0$, the pre-image of $\gamma _{\rho }$ passing through $%
s $ must contain two different components $\eta _{k,j},$ we conclude that
the points $v_{k,j}$ are the only zeros of $\zeta ^{\prime }$ having images
by $\zeta $ situated in the unit disc$.$ On the other hand a zero of $\zeta
^{\prime }$ in some strip $S_{k}$ with image by $\zeta $ outside the unit
disc cannot belong to components of pre-images of $\gamma _{\rho }$ with
different values of $\rho ,$ due to the fact that $\zeta $ is a single
valued function$.$ The only way for a branch point of $\zeta $ to have the
image on a $\gamma _{\rho }$ with $\rho >1$ is when a component of the
pre-image of $\gamma _{\rho }$ turns around the bounded component of the
pre-image of the unit circle and has a self intersection point.

\bigskip

\begin{center}
\includegraphics[width=1.8truein,height=2truein]{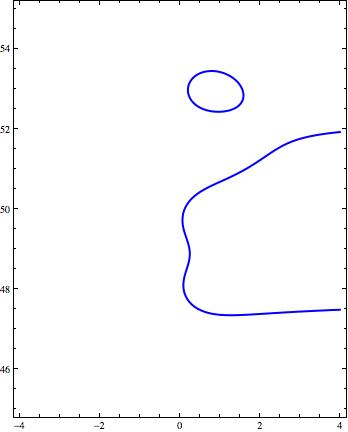} \ \ \ %
\includegraphics[width=1.8truein,height=2truein]{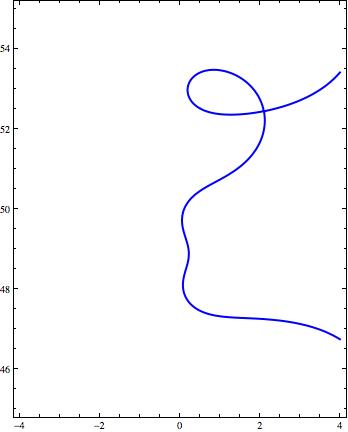}\ \ \ %
\includegraphics[width=1.8truein,height=2truein]{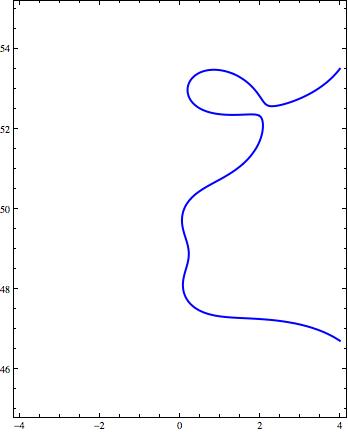}\\[0pt]
\ \\[0pt]
(a)\hspace{1.8in} (b)\hspace{1.8in} (c) \\[0pt]
\ \\[0pt]
\includegraphics[width=2truein,height=2truein]{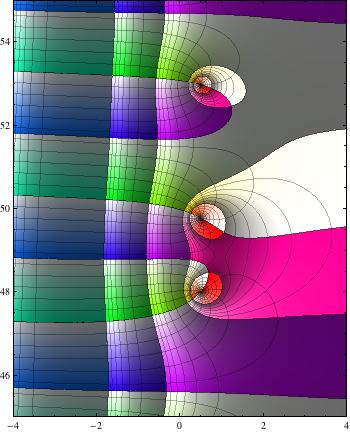}\ \ \ \ \ \ \ \ %
\includegraphics[width=1.8truein,height=2truein]{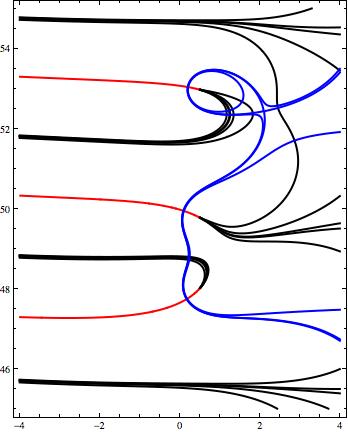}\ \ \\[0pt]
(d)\hspace{2in} (e)\\[0pt]
\ \\[0pt]
Figure 9
\end{center}

Fig. 9 illustrates the situation where a component of the pre-image of $%
\gamma _{\rho }$ has a self-intersection point. In the box [-4,4] x [45,55],
Fig. 9(a), two components of the pre-image of the circle $\gamma _{\rho }$
with $\rho =1$ are visible: a bounded one on the upper part of the box,
containing a unique non trivial zero of $\zeta $ and an unbounded one
covering the right lower corner of the box and containing two non trivial
zeros. As the radius $\rho $ takes values greater than $1$, the two
components expand, Fig. 9(b,c), touching each other for $\rho =\rho
_{0}\approx 1.042$ in Fig. 9(b). We can interpret the pre-image of $\gamma
_{\rho _{0}}$ as having a unique unbounded component with a
self-intersection point $v_{5,2}.$ It borders three domains, one bounded and
two unbounded. As $\rho $ takes values greater than $\rho _{0},$ the bounded
component \emph{opens}, Fig. 9(c) and we get a unique unbounded component
bordering two unbounded domains. It is obvious that $v_{5,2}$ is a branch
point of $\zeta .$ Indeed, the arcs of the pre-image of $\gamma _{\rho _{0}}$
situated in a small neighborhood of $v_{5,2}$ are mapped by $\zeta $ onto an
arc of $\gamma _{\rho _{0}}$ containing $\zeta (v_{5,2})$. Thus $\zeta
^{\prime }(v_{5,2})=0$ and $v_{5,2}$ is a simple zero of $\zeta ^{\prime }.$
Fig. 9(e) is a superposition of Fig. 6 and Fig. 9(a-c) showing the domains
mapped by $\zeta $ outside the circles $\gamma _{\rho }$ and the sectors in
Fig. 6. It helps locate $v_{5,2}$ on Fig. 9(d).\bigskip

It is obvious that the scenario described in Fig. 9 happens for every
bounded component of the pre-image of the unit circle. In other words, to
every $u_{k,j},$ except for $u_{k,0},$ corresponds a branch point $v_{k,j}$
of $\zeta $ situated in the strip $S_{k}.$ We notice that components of the
pre-image of $\gamma _{\rho }$ with different values of $\rho $ cannot
intersect, since this would contradict the single value nature of $\zeta .$
Thus, $\zeta $ cannot have branch points other than $v_{k,j}.$ This remark
allows us to state:\bigskip

\textbf{Theorem 5}. \textit{All the non real zeros of }$\zeta ^{\prime }$%
\textit{\ are situated in the right half plane.}\bigskip

Proof: Indeed, those zeros are the points $v_{k,j}$ previously obtained, and
their location in the right half plane is obvious.\bigskip\ There can be no
other non real zero.

Let us connect $u_{k,j}$ and $v_{k,j}$ by disjoint Jordan arcs $%
L_{k,j}^{\prime }$ . Then $\gamma _{k,j}=\zeta (L_{k,j}^{\prime })$ are
Jordan arcs connecting $z=1$ with some points $z=\zeta (v_{k,j}).$ We
perform continuations over $\gamma _{k,j}$ from every $v_{k,j}$ in the
opposite direction of $L_{k,j}^{\prime }$ obtaining unbounded curves $%
L_{k,j}^{\prime \prime }$ such that $\lim_{\sigma \rightarrow +\infty }\zeta
(\sigma +it)=1,$ as $\sigma +it\in L_{k,j}^{\prime \prime }.$ Let us denote $%
L_{k,j}=L_{k,j}^{\prime }\cup L_{k,j}^{\prime \prime }.$ When a point $s$
travels on $L_{k,j}$ from $u_{k,j}$ to $\infty ,$ the point $\zeta (s)$
travels on $\gamma _{k,j}$ from $z=1$ to $z_{k,j}=\zeta (v_{k,j})$ and back
to $z=1.$Two curves $L_{k,j}$ and $L_{k,j^{\prime }}$cannot meet each other
since the corresponding arcs $L_{k,j}^{\prime }$ and $L_{k,j^{\prime
}}^{\prime }$ are supposed to be disjoint. If we add to every $L_{k,j}$ the
part of $\Gamma _{k,j}$ corresponding to the interval $(1,+\infty )$ we
obtain curves $C_{k,j}$ which are unbounded in both directions.\ Obviously,
these curves do not intersect and they
divide every $m$-strip $S_{k}$ into $m+1$ strips $\Omega _{k,j}$, $%
j=0,1,2,\ldots ,m.$ Let us show that $\Omega _{k,j}$ are fundamental domains
of $\zeta .$\bigskip

Four types of domains $\Omega _{k,j},$ depending on their boundaries, can be
distinguished. These boundaries can be: $\Gamma _{1}^{\prime }$ and $\Gamma
_{2}^{\prime }$, and for $k\geq 2,$ a $\Gamma _{k}^{\prime }$ and the
corresponding $C_{k,1},$ a $C_{k,j}$ and $C_{k,j+1},$ for $1\leq j<m$ and
finally $C_{k,m}$ and $\Gamma _{k+1}^{\prime }.$ We notice that every domain 
$\Omega _{k,j}$ is mapped conformally by $\zeta $ onto the complex plane
with a slit, thus it is indeed a fundamental domain. The slit is the
interval $[1,+\infty )$ for the first domain, $[1,+\infty )\cup \gamma
_{k,1},$ respectively $[1,+\infty )\cup \gamma _{k,m}$ for the second and the
last type of domains and finally $[1,+\infty )\cup \gamma _{k,j}\cup \gamma
_{k,j+1}$ for the third type of domains.

\bigskip

\begin{center}
\includegraphics[width=1.8truein,height=2.2truein]{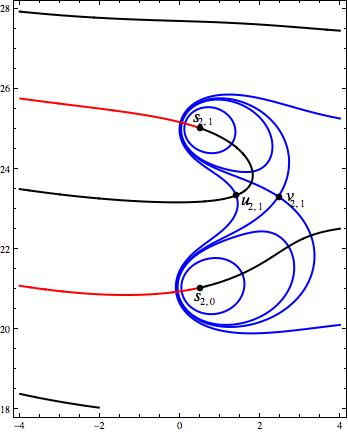} \ \ \ %
\includegraphics[width=1.8truein,height=2.2truein]{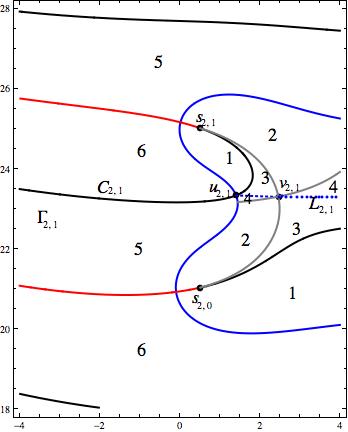}\ \ \ %
\includegraphics[width=1.8truein,height=1.2truein]{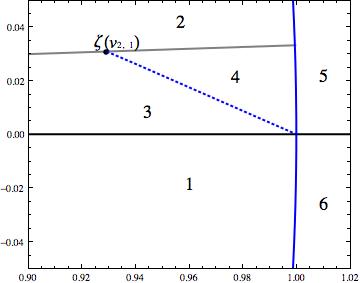}\\[0pt]
\hspace{0in}(a)\hspace{1.8in}(b)\hspace{1.8in}(c)\\[0pt]
\ \\[0pt]
Figure 10
\end{center}

\bigskip

Fig. 10 describes the conformal mapping by $\zeta $ in the strip $S_{2}.$
Fig. 10(a) shows the way in which $v_{2,1}$ is obtained as the point where
the components of the pre-image of the circle $\gamma _{\rho }$ meet each
other for $\rho =\rho _{0}\approx 0.9296.$ In Fig. 10(b) a curve $C_{2,1}$
is shown dividing the strip $S_{2}$ into two fundamental domains. Details of
the conformal mapping by $\zeta $ of these domains onto the complex plane
with a slit are visible when one compares the Fig. 10(b) and the Fig.
10(c).\bigskip

The question has been asked whether some of the non trivial zeros of $\zeta $
can be multiple zeros and a few papers have been published on this topic. In
our knowledge there are just statistical estimations of the proportion of
such zeros [7]. The method of fundamental domains allows us to give a
definite answer to this question. Indeed, every non trivial zero belongs to
a fundamental domain $\Omega _{k,j},$ which is mapped bijectively by $\zeta $
onto the complex plane with a slit. Therefore, we have\bigskip :

\textbf{Theorem 6}.\textit{\ All the zeros of }$\zeta $\textit{\ are simple
zeros.\bigskip\ }

Having in view the remark after Theorem 4, we can restate that theorem by
saying that all the zeros of $\zeta ^{\prime }$ are simple
zeros.\bigskip

\bigskip \textbf{7. The Group of Covering Transformations of }$(\mathbb{C}%
,\zeta )$\bigskip

In order to find the group of covering transformations of $(\mathbb{C},\zeta
)$ we need to rename the fundamental domains. We proceed in a way similar to
that of Gamma function. Let us denote by $\sigma _{n}$ the branch points of $%
\zeta $ situated on the negative real half axis counted in an increasing
order of their module. Let\ $\Omega _{0}$ be the domain bounded by the
interval $(\sigma _{2},+\infty ),$ the branch from the upper half plane of
the component passing through $\sigma _{2}$ of the pre-image of the negative
real half axis and $\Gamma _{1}^{\prime }$. We notice that $\Omega _{0}$ is
mapped conformally by $\zeta $ onto the complex plane with a slit alongside
the real axis complementary to the interval $(\zeta (\sigma _{1}),1),$ where 
$\sigma _{1}$ is the first branch point of $\zeta .$ The domains $\Omega
_{-n},$ $n\in N$ are the regions from the upper half plane bounded by the
real axis and the branches of the pre-image of the negative half axis
starting at $\sigma _{2n},$ respectively $\sigma _{2(n+1)}$. They are mapped
conformally by $\zeta $ onto the complex plane with a slit alongside real
axis from $-\infty $ to $\zeta (\sigma _{2n+1}).$ Finally the domains $%
\Omega _{n},$ $n\in N$ are the former domains $\Omega _{k,j}$ counted
starting from the positive real half axis and going up. All these domains
are fundamental domains for $\zeta .$ So are the domains $\widetilde{\Omega }%
_{n}$ symmetric to them with respect to the real axis. \bigskip

We define, as in the case of the Gamma function mappings $U_{k}$ and and $H$
\ by the formulas similar to $(6)$ \ and $(8),$ where $\Gamma $ is replaced
by $\zeta $ and notice that the group generated by $U_{1}$ \ and $H$ is the
group $G$ of covering transformations of $(\mathbb{C},\zeta ).$ The group
generated by $U_{1}$ is an infinite cyclic subgroup of $G.$\bigskip

\textbf{Acknowledgments: }The authors are grateful to Cristina Ballantine
for providing computer generated graphics.

\textbf{\ }

\textbf{References\bigskip }

[1] Ahlfors, L.V., Complex Analysis, International Series in Pure and
Applied Mathematics (1979)

[2] Andreian Cazacu, C. and Ghisa, D., Global Mapping Properties of Analytic
Functions, Proceedings of the 7-th ISAAC Congress, London, U.K. (2009)

[3] Ballantine, C. and Ghisa, D., Color Visualization of Blaschke Product
Mappings, to appear in Complex Variables and Elliptic Equations

[4] Ballantine, C. and Ghisa, D., Color Visualization of Blaschke
Self-Mappings of the Real Projective Plan, to appear in Proceedings of the
International Conference on Complex Analysis and Related Topics, Alba Iulia,
Romania (2008)

[5] Ballantine, C. and Ghisa, D., Global Mapping Properties of Rational
Functions, to appear in Proceedings of the 7-th ISAAC Congress, London,
U.K., (2009)

[6] Barza, I. and Ghisa, D., The Geometry of Blaschke Product Mappings,
Proceedings of the 6-th International ISAAC Congress, Ankara (2007) in
Further Progress in Analysis, Editors H.G.W.Begehr, A.O. Celebi \&
R.P.Gilbert, World Scientific, 197--207 (2008)

[7] \ Cheer, A.Y. and Goldston, D.A., Simple zeros of Riemann Zeta function,
Proc. of the AMS, Vol 118, No.2, 365-372 (1993)

[8] http://en.wikipedia.org/wiki/Gamma\_function

\end{document}